\documentclass[11pt]{article}
\usepackage{amsfonts}
\usepackage{amssymb}
\usepackage{amsmath}
\def\sqr#1#2{{\vcenter{\vbox{\hrule height.#2pt
              \hbox{\vrule width.#2pt height#1pt \kern#1pt \vrule
width.#2pt}
              \hrule height.#2pt}}}}
\def\dbE{{\mathbb{E}}}
\def\dbF{{\mathbb{F}}}

\def\dbN{{\mathbb{N}}}

\def\dbP{{\mathbb{P}}}

\def\dbR{{\mathbb{R}}}

\def\3n{\negthinspace \negthinspace \negthinspace }
\def\2n{\negthinspace \negthinspace }
\def\1n{\negthinspace }

\def\ds{\displaystyle}
\def\O{\Omega}

\def\cF{{\cal F}}

\def\ms{\medskip}

\def\cd{\cdot}

\def\({\Big (}
\def\){\Big )}
\def\[{\Big[}
\def\]{\Big]}

\def\={\buildrel \triangle \over =}

\def\be{\begin{equation}}
\def\bel{\begin{equation}\label}
\def\ee{\end{equation}}
\def\bea{\begin{eqnarray}}
\def\eea{\end{eqnarray}}
\def\bt{\begin{theorem}}
\def\et{\end{theorem}}
\def\bc{\begin{corollary}}
\def\ec{\end{corollary}}
\def\bl{\begin{lemma}}
\def\el{\end{lemma}}
\def\bp{\begin{proposition}}
\def\ep{\end{proposition}}
\def\br{\begin{remark}}
\def\er{\end{remark}}
\def\ba{\begin{array}}
\def\ea{\end{array}}
\def\bd{\begin{definition}}
\def\ed{\end{definition}}

\newtheorem{lemma}{Lemma}[section]
\newtheorem{remark}{Remark}[section]

\newtheorem{theorem}{Theorem}[section]
\newtheorem{corollary}{Corollary}[section]

\newtheorem{definition}{Definition}[section]
\newtheorem{proposition}{Proposition}[section]

\makeatletter
   
   \@addtoreset{equation}{section}
\makeatother

\begin{document}

\title{\bf Numerical Solutions of Backward Stochastic Differential Equations: A Finite Transposition Method}
\author{Penghui Wang\thanks{School of Mathematics, Shandong University, Jinan, 250100, China. {\small\it e-mail:} {\small\tt phwang@sdu.edu.cn}. \ms}
 ~~~
  and~~~
Xu Zhang\thanks{Key Laboratory of Systems and Control, Academy of
Mathematics and Systems Science, Chinese Academy of Sciences,
Beijing 100190, China; Yangtze Center of Mathematics, Sichuan
University, Chengdu 610064, China. {\small\it e-mail:} {\small\tt
xuzhang@amss.ac.cn}.}}
\maketitle
\begin{abstract}

In this note, we present a new numerical method for solving backward
stochastic differential equations. Our method can be viewed as an
analogue of the classical finite element method solving
deterministic partial differential equations.

\end{abstract}

\section{Introduction}

Linear and nonlinear Backward Stochastic Differential Equations
(BSDEs in short) were introduced in \cite{Bis1} and \cite{PP},
respectively. It is well-known that BSDE plays crucial roles in
Stochastic Control, Mathematical Finance etc. Clearly, for
applications, it deserves to develop effective numerical methods for
BSDEs.

Let $T>0$ and   $(\O,\cF,\dbF,\dbP)$ be a complete filtered
probability space with $\dbF=\{\cF_t\}_{t\in[0,T]}$, on which a
$1$-dimensional standard Brownian motion $\{w(t)\}_{t\in[0,T]}$ is
defined. We denote by $L_{\cF_t}^2(\O;\dbR^n)$  ($n\in \dbN$) the
Hilbert space consisting of all $\cF_t$-measurable ($\dbR^n$-valued)
square integrable random variables; by
$L^2_{\dbF}(\Omega;L^r(0,T;\dbR^n))$ ($1\leq r\leq \infty$) the
Banach space consisting of all $\dbR^n$-valued $\{\cF_t\}$-adapted
processes $X(\cdot)$ such that
$\dbE|X(\cdot)|_{L^r(0,T;\dbR^n)}^2<\infty$; and by
$L^{2}_{\dbF}(\O;D([0,T];\dbR^n))$ the Banach space consisting of
all $\dbR^n$-valued $\{ {\cF}_t \}$-adapted c\`adl\`ag processes
$X(\cdot)$ such that
$\mathbb{E}(|X(\cdot)|^2_{L^{\infty}_{\dbF}(0,T;\dbR^n)}) < \infty$.
For $y_T\in L^2_{\cF_T}(\Omega;\dbR^n)$ and $f(\cdot,\cdot,\cdot)$
satisfies $f(\cdot,0,0)\in L^2(\Omega;L^1(0,T;\dbR^n))$ and the
usual globally Lipschitz condition, we consider the following BSDE
\begin{eqnarray}\label{eq1.1}
\left\{
\begin{array}{l}
dy(t)=f(t,y(t),Y(t))dt+Y(t)dw(t)\quad in\ [0,T],\\
y(T)=y_T.
\end{array}
\right.
\end{eqnarray}
Various numerical methods have been developed to solve equation
(\ref{eq1.1}), say in \cite{BT,MPST,PX,Zhang} and the references
therein. These methods use essentially the strong form of
(\ref{eq1.1}), which holds true only if $\dbF$ is the natural
filtration generated by the Brownian motion. Also, it seems that the
previous methods need to compute the conditional expectation, which
is in general not easy to be furnished numerically.

In this Note, we shall present a new numerical method solving BSDEs
from the viewpoint of transposition solution introduced in
\cite{LZ}, as recalled below.\vskip1.5mm

\begin{definition}\label{OK1}
A couple $(y(\cdot), Y(\cdot))\in
L^2_{\dbF}(\Omega;D([0,T];\dbR^n))\times
L^2_{\dbF}(\Omega;L^2(0,T;\dbR^n))$ is called a transposition
solution of BSDE $($\ref{eq1.1}$)$, if for any $t\in[0,T]$,
$(u(\cdot),v(\cdot),\eta)\in L^2_{\dbF}(\Omega;
L^1(t,T;\dbR^n))\times L^2_{\dbF}(\Omega; $ $L^2(t,T;\dbR^n))\times
L^2_{\cF_t}(\Omega;\dbR^n)$, the following variational equation
holds
\begin{eqnarray}\label{eq1.2}
&&\dbE\langle z(T),y_T\rangle -\dbE\langle \eta,y(t)\rangle\nonumber\\
&&=\dbE\int_t^T\langle z(\tau),
f(\tau,y(\tau),Y(\tau))d\tau+\dbE\int_t^T \langle
u(\tau),y(\tau)\rangle d\tau+\dbE\int_t^T\langle v(\tau),
Y(\tau)\rangle d\tau,
\end{eqnarray}
where $z(\tau)=\eta+\int_t^\tau u(s)ds+\int_t^\tau v(s)dw(s)$.
\end{definition}
\vskip1.5mm

We refer to \cite{LZ} for the well-posedness of equation
$($\ref{eq1.1}$)$ in the sense of transposition solution. It is easy
to see that, if this equation admits a strong solution (say under
the assumption of natural filtration), it coincides with the
transposition solution.

Our numerical schemes for solving equation $($\ref{eq1.1}$)$ can be
described as follows. \vskip1.5mm
\begin{itemize}
\item[1)] Take a suitable finite dimensional subspace $H_m$ of $L^2_{\dbF}(\Omega;L^2(0,T;\dbR^n))$;

\item[2)] If the solution of $($\ref{eq1.1}$)$ exists, then it should satisfy the variational equation (\ref{eq1.2}) for any $u,v\in H_m$.
By taking $u$ and $v$ to be the orthonormal basis of $H_m$, we
obtain a system of approximating equations;
\item[3)] If the solution of the  system of approximating equations exists, then we find a class of numerical solutions of
$($\ref{eq1.1}$)$;

\item[4)] Finally, we show the convergence of the above numerical solutions.
\end{itemize}
\vskip1.5mm

Clearly, the above procedure is, in spirit, very close to that of
the classical finite element method solving deterministic PDEs
(e.g., \cite{BR}). Therefore, our method  to solve BSDEs can be
viewed as a stochastic version of the finite element-type method.
Nevertheless, the notion of ``stochastic finite element method" has
already been used for other purpose, say \cite{GS,KH,N} and
references therein for solving random PDEs. Note also that our
method is quite different from that in these references, and
therefore instead we call it a finite transposition method.

There are at least two reasons for us to develop this new numerical
approach for BSDE. The first one is that, we can solve the BSDE with
general filtration. The second is that, in our approach,  we do not
need to compute the conditional expectation.

We refer to \cite{WZ} for the details of the proofs of the results
in this Note and other results in this context.

\section{Numerical schemes and convergence}
For simplicity, we consider only the following linear BSDE (with
$f(\cdot)\in L^2_{\dbF}(\Omega;L^1(0,T;\dbR^n))$)
\begin{eqnarray}\label{bs1}
\left
\{\begin{array}{l}
dy(t)=f(t)dt+Y(t)dw(t),\ \ t\in [0,T),\\
y(T)=y_T.
\end{array}
\right.
\end{eqnarray}

Assume that $L_{\cF_T}^2(\Omega;\dbR^n)$ is a separable Hilbert
space. For any $N\in\dbN$, write $ \mathfrak{R}_{N}=
\{t_\ell\,|\,t_\ell={\ell\over 2^N}T,\ \ell=0,\cdots, 2^N\}$. For
any $k\in\{0,\cdots, 2^N-1\}$, define a sequence of simple processes
$\{e_{ki}(\cdot,\cdot)\}_{i=1}^{M_{k,N}}$ by
\begin{eqnarray}
e_{ki}(t,\omega)=\left\{\begin{array}{ll} \chi_{[t_k,t_{k+1})}(t)h_{ki}(\omega), & 0\leq k<2^N-1,\\  \chi_{[t_{k},T]}(t)h_{ki}(\omega), & k=2^N-1,\end{array}\right.
\end{eqnarray}
where $ \{M_{0,N},M_{1,N},\cdots,M_{2^N-1,N}\}$ is an increasing
sequence of integers. Since $L^2_{\cF_t}(\Omega;\dbR^n)\subseteq
L^2_{\cF_s}(\Omega;\dbR^n)$ (for any $0\le t<s\le T$), we may assume
that $\{h_{ki}\}$ satisfy the following:

\begin{itemize}
\item[1)] For any fixed $k\in\{0,\cdots, 2^N-1\}$, $\{h_{ki}\}_{i=0}^{M_{k,N}}$ is an orthogonal set in $L^2_{\cF_{t_k}}(\Omega;\dbR^n)$,
and the norm
$|h_{ki}|_{L^2_{\cF_{t_k}}(\Omega;\dbR^n)}=\sqrt{2^N\over T}$, and
hence $|e_{ki}|_{L^2_{\dbF}(\Omega;L^2(0,T;\dbR^n))}=1$;

\item[2)] If $0\le k<\ell\leq 2^N-1$, then
$\{h_{ki}\}_{i=0}^{M_{k,N}}\subset\{h_{\ell i}\}_{i=0}^{M_{\ell,N}}$; and

\item[3)] If $s_0={k_0\over 2^{N_0}}T$ for some $N_0\in \mathbb N$ and $k_0\in\{0,\cdots, 2^{N_0}-1\}$, then $s_0=2^{\ell-N_0}k_0 T/2^\ell\in \mathfrak{R}_\ell$ for any $\ell\geq N_0$.
For $\ell\geq N_0$, write $k_\ell=2^{\ell-N_0}k_0$. Then,
$\{h_{k_ji}\}_{i=0}^{M_{k_j,j}}\subset
\{h_{k_\ell i}\}_{i=0}^{M_{k_\ell,\ell}}$ for $N_0\leq j<\ell$, and
$\bigcup\limits_{\ell=N_0}^\infty\{h_{k_\ell i}\}_{i=1}^{M_{k_\ell,\ell}}$ is an
orthogonal basis of $L^2_{\cF_{s_0}}(\Omega;\dbR^n)$.
\end{itemize}

Denote by $H_N$ the subspace of $L^2_{\dbF}(\Omega;L^2(0,T;\dbR^n))$
spanned by
$\{e_{0i}\}_{i=1}^{M_{0,N}},\cdots,\{e_{2^N-1,i}\}_{i=1}^{M_{2^N-1,N}}$.
This is the finite element subspace that we will employ below.
Replace $L^2_{\dbF}(\Omega;L^2(0,T;\dbR^n)$ by $H_N$, then the
desired numerical scheme follows by trying to find  $y_N,Y_N\in H_N$
such that (\ref{eq1.2}) (with $f(\cd,y(\cd),Y(\cd))$ replacing by
$f(\cd)$) holds for $\eta=0$ and for all $ u,v\in H_N$.

To find the $y_N$ in $H_N$, suppose
$y_N=\sum\limits_{k=0}^{2^N-1}\sum\limits_{i=0}^{M_{k,N}}\alpha_{ki}e_{ki}$.
Choosing $u=e_{ki}$, $v=0$ and $\eta=0$, we get $z_{ki}(t)=\int_0^t
e_{ki}(\tau) d\tau$, and hence $ \dbE\langle z_{ki}(T),
y_T\rangle=\dbE\int_0^T\langle z_{ki}(\tau), f(\tau)\rangle
d\tau+\sum\limits_{\ell,j} \alpha_{\ell j}\dbE\int_0^T\langle
e_{ki}(\tau), e_{\ell j}(\tau)\rangle d\tau. $ Since $\{e_{ki}\}$ is
an orthonormal basis of $H_N$, it follows that $ \dbE\int_0^T\langle
e_{\ell j}(\tau),e_{ki}(\tau)\rangle d\tau=\delta_{k\ell}\delta_{ij}. $
Therefore,
\begin{eqnarray}
\alpha_{ki}=\frac{T}{2^N}\dbE\langle h_{ki},
y_T\rangle-\dbE\int_0^T\langle (\tau\wedge t_{k+1}-\tau\wedge
t_k)h_{ki}, f(\tau)\rangle d\tau.
\end{eqnarray}
Similarly, suppose
$Y_N=\sum\limits_{k=0}^{2^N-1}\sum\limits_{i=0}^{M_{k,N}}\beta_{ki}e_{ki}$.
By taking $u=0,\eta=0$ and $v=e_{ki}$ to get a corresponding
$z_{ki}(t)=\int_0^t e_{ki}(\tau)dw(\tau)$, we find that
\begin{eqnarray*}
\beta_{ki}=\dbE\langle (w(t_{k+1})-w(t_{k}))h_{ki}, y_T\rangle - \dbE\int_0^T\langle (w(\tau\wedge t_{k+1})-w(\tau\wedge t_k))h_{ki},f(\tau)\rangle d\tau.
\end{eqnarray*}

We now show the convergence of the sequence $\{(y_N, Y_N)\}$ of
numerical solutions constructed above.\vskip1.5mm
\begin{theorem}
Let $(y, Y)$ be the transposition solution of (\ref{bs1}). Then
$y_N$ and $Y_N$ are projections of $y$ (viewing as an element of
$L^2_{\dbF}(\Omega;L^{2}(0,T;\dbR^n))$) and $Y$ onto $H_N$,
respectively.
\end{theorem}\vskip1.5mm

By the construction of $H_N$, it is clear that $H_N\subset H_M$
provided $N<M$.  Moreover, $\ds\bigcup_{N=1}^\infty H_N$ is dense in
$L^2_{\dbF}(\Omega;L^{2}(0,T;\dbR^n))$. Hence,
\begin{eqnarray}\label{eq1}
\dbE\int_0^T |y_N(\tau)-y(\tau)|^2 d\tau+\dbE\int_0^T
|Y_N(\tau)-Y(\tau)|^2 d\tau\to 0,\  as\ N\to\infty.
\end{eqnarray}
Furthermore, starting from (\ref{eq1}), we can show the following
convergent result.\vskip1.5mm
\begin{theorem}
As $ N\to\infty$, $(y_N, Y_N)$ tends to $(y, Y)$ in
$L^2_{\dbF}(\Omega;D([0,T],\dbR^n))\times
L^2_{\dbF}(\Omega;L^{2}(0,T;\dbR^n))$.
\end{theorem}
\section*{Acknowledgements}
This work is partially supported by the NSF of China under grants
10901032 and 10831007,  the National Basic Research Program of China
(973 Program) under grant 2011CB808002, and by Innovation Foundation
of Shandong University. The second author thanks Dr. Qi Zhang (Fudan
University) for stimulating discussions.

\end{document}